\documentclass[12pt,a4paper]{article}

\usepackage{graphicx}
\usepackage{latexsym}
\usepackage{amsfonts}

\author{Christophe Bavard, K\'aroly J. B\"or\"oczky\footnote{Supported by
OTKA grants 068398 and 75016, and by the EU Marie Curie TOK
project DiscConvGeo and FP7 IEF grant GEOSUMSETS},\\
Borb\'ala Farkas, Istv\'an Prok, Lluis Vena, Gergely Wintsche}
\title{Equality
in L\'aszl\'o Fejes T\'oth's triangle bound for hyperbolic surfaces}

\newcommand{\HH}{\mathbb{H}}
\newcommand{\R}{\mathbb{R}}

\begin{document}

\maketitle

\newtheorem{theo}{Theorem}[section]
\newtheorem{coro}[theo]{Corollary}
\newtheorem{lemma}[theo]{Lemma}
\newtheorem{remark}[theo]{Remark}
\newtheorem{prop}[theo]{Proposition}
\newtheorem{conj}[theo]{Conjecture}
\newtheorem{example}[theo]{Example}

\begin{abstract}
For $k\geq 7$, we determine the minimal area of
a compact hyperbolic surface, and an
oriented compact hyperbolic surface that can be tiled by embedded regular
triangles of angle $2\pi/k$. Based on this, all the cases of equality
in  L\'aszl\'o Fejes T\'oth's triangle bound for hyperbolic surfaces are described.
\end{abstract}

\section{Introduction}
\label{secmain}

Let $\HH^2$ be the hyperbolic plane of curvature $-1$
(see J.G. Ratcliffe \cite{Rat94} for
facts and references on hyperbolic geometry,
and D.V. Alekseeskij, E.B. Vinberg, A.S. Solodovnikov \cite{AVS93}
for an in depth study on discrete groups of isometries).
A compact hyperbolic surface  $X$ can be
obtained as the quotient of $\HH^2$ by a discrete
group $\Gamma$ of isometries that acts fixed point free
on $\HH^2$. For a regular polygon $\Pi$ in $\HH^2$,
we say that $X$ can be tiled
by copies of $\Pi$ if there exists an edge to edge
tiling of $\HH^2$ by congruent
copies of $\Pi$ that is invariant under $\Gamma$.
The equivalence classes of the
tiles, edges, and vertices of the
tiling in $\HH^2$ with respect to $\Gamma$
correspond to the tiles, edges, and vertices,
respectively, of the
tiling on $X$ (as a CW-decomposition
of $X$). In this
case there exists an integer $p\geq 3$
such that the angles of $\Pi$ are $2\pi/p$. If $\Pi$
is a regular triangle then $p\geq 7$.

In this note, we provide a simple proof
for the following statement.

\begin{theo}
\label{main}
For $k\geq 7$, the minimal area of
a compact hyperbolic surface that can be tiled by embedded regular
triangles of angle $2\pi/k$, $k\geq 7$, is $N(k-6)\frac{\pi}3$,
where $N$ is the minimal positive integer
such that $Nk$ is divisible by six. If the
surface is assumed to be orientable then the
minimal area is the same if $k\equiv 2,6,10\,{\rm mod}\, 12$, and
twice the previous value otherwise.
\end{theo}
{\bf Remark } It follows by
A.L. Edmonds, J.H. Ewing, R.S. Kulkarni \cite{EEK82i}
that there is a compact hyperbolic surface that can be 
tiled by $t$ embedded regular
triangles of angle $2\pi/k$, $k\geq 7$, if
and only if $t$ is even, and $3t$ is divisible by $k$.\\ 

We note that Theorem~\ref{main} is proved in
a more general setup, and using a more
involved construction in 
A.L. Edmonds, J.H. Ewing, R.S. Kulkarni \cite{EEK82a}.
In addition, the case of Theorem~\ref{main}
when $k$ is
divisible by six, is proved by C. Bavard
\cite{Bav96}.

In other words, for $k\geq 7$, let $\Gamma_k$ be the
isometry group of the edge to edge tiling
by regular triangles of angle $2\pi/7$. Then
Theorem~\ref{main} yields (using its notation)
that the minimal index
of a subgroup of $\Gamma_k$ that
acts fixed point free on $\HH^2$ is $2Nk$.

Our motivation to consider Theorem~\ref{main} 
comes from the theory of packing  and
 covering by equal circular discs on
a compact surface $X$ of constant curvature.
Let $T$ be a regular triangle of angle $\alpha$
in the universal covering surface, which is
$\HH^2$, $\R^2$, and the unit sphere $S^2$
in $\R^3$ if the constant curvature is $-1$, $0$
or $1$, respectively. In addition, let 
$r$ and $R$ be inradius and circumradius,
 respectively, of $T$, and let 
$A(\cdot)$ denote area.
According to the triangle bound due to
 L. Fejes T\'oth \cite{FTL64} (see also
G. Fejes T\'oth, W. Kuperberg \cite{FTK93} for a survey), if $X$ can be
 covered by $n$ metric balls of radius $R$,
 then
 $$
 A(X)\leq n\cdot\frac{2\pi}{3\alpha}\cdot A(T).
 $$
In addition, if $X$ contains $n$ pairwise disjoint
open embedded balls of radius $r$, then
$$
 A(X)\geq n\cdot\frac{2\pi}{3\alpha}\cdot A(T).
 $$
 If the injectivity radius of $X$ is at least $r$,
 then $A(X)\geq \frac{2\pi}{3\alpha}\cdot A(T)$.
In these bounds, equality holds if and
only if $\alpha=2\pi/k$ for some integer $k\geq 3$,
and $X$ can be tiled by copies of $T$
using exactly $n$ vertices.
 To characterize the cases when
we have equality in the estimate
using the injectivity radius ($n=1$)
is the subject of C. Bavard \cite{Bav96}.

To charaterize regular triangulations when the curvature of $X$ is non-negative,
is rather straightforward.
If the curvature of $X$ is positive, then
$k\in\{3,4,5\}$, as the sum of the angles of $T$
is larger, than $\pi$. There are two types of possible $X$.
If $X$ is a sphere then the possible regular triangulations are
determined by the platonic solids: tetrahedron, octahedron and
icosahhedron.  If $X$ is a projective plane then 
the possible regular triangulations are obtained by
identifying the opposite faces of the octahedron and the 
icosahhedron.

Next, there are exactly two types of surfaces
of curvature zero; namely, tori and Klein bottles.
In these cases $k=6$.
For any $n\geq 1$, it is easy to construct a
surface of either type with a regular
triangulation having $n$ vertices.

Therefore, it remains to find regularly triangulated hyperbolic surfaces,
when $k\geq 7$. This is the content
of   Theorem~\ref{main} and the following remark.

\section{Some basic properties of hyperbolic surfaces}

First, we provide some properties of
tilings of a compact hyperbolic surface $X$.
The Euler characteristic of $X$ is denoted by $\chi(X)$,
which is negative, and is always even in the case of
oriented surfaces.
It follows from the Gau{\ss}-Bonnet formula that
the area of $X$ is $2\pi|\chi(X)|$, therefore
minimizing the area is equivalent to
maximizing the Euler characteristic.

\begin{lemma}
\label{nk6}
If a compact hyperbolic surface  $X$
is tiled by embedded regular
triangles of angle $2\pi/k$, $k\geq 7$, and the
tiling has $n$ vertices (as a CW-decomposition), then
 $nk$ is divisible by six, and
$$
\chi(X)=n\cdot\left(1-\frac{k}6\right).
$$
\end{lemma}
{\it Proof: } Since each vertex
of the tiling is incident with exactly $k$
edges, the tiling has $nk/3$ faces
and $nk/2$ edges. Therefore, we have
$$
\chi(X)=n-\frac{nk}2+\frac{nk}3=
n\cdot\left(1-\frac{k}6\right).
$$
Since $nk$ is divisible
by two and three, it is divisible
by six. Q.E.D.

After the necessary condition above,
we present some ways to construct
suitable hyperbolic surfaces.
Using the oriented double cover
of a non-orientable surface, we have
the following.

\begin{lemma}
\label{nonoriented}
If a compact non-orientable hyperbolic surface  $X$
is tiled by embedded regular
triangles of angle $2\pi/k$, $k\geq 7$, using
$n$ vertices, then
there exists a compact oriented hyperbolic surface
$\tilde{X}$, which
is tiled by embedded regular
triangles of angle $2\pi/k$ using
$2n$ vertices, and
$\chi(\tilde{X})=2\chi(X)$.
\end{lemma}

Assume that a compact hyperbolic surface  $X$
is tiled by embedded regular
triangles of angle $2\pi/k$, $k\geq 7$, and 
the tiling has $n$ vertices.
By duality, it is equivalent that
$X$ can be tiled by $n$ embedded
regular $k$-gons of angle $2\pi/3$.

To obtain $X$, one may start with $n$
regular $k$-gons, and identify the suitable
pairs out of the total
$nk$ edges of the $n$-gons.
We describe a combinatorial way to encode
pairings that we need in the sequel.

Let $\Pi_1,\ldots,\Pi_n$ be convex
polygons with at least seven vertices such that the their
total number of vertices (or edges)
is $3m$ for some even integer $m$.
A {\it proper labeling} of
the vertices of  $\Pi_1,\ldots,\Pi_n$
means an assignment of a label
from $\{1,\ldots,m\}$ to each vertex
with the following properties.

\begin{description}
\item[(i)] Any  $i\in\{1,\ldots,m\}$
is the label of exactly three of the vertices.
\item[(ii)] For any non-empty proper
subset of $\Pi_1,\ldots,\Pi_n$ there
exists a label $i\in\{1,\ldots,m\}$ that occurs
once or twice as the label of some vertex of
a convex polygon in the subset.
\item[(iii)] If a $\Pi_l$ has
two consecutive vertices with label $i\in\{1,\ldots,m\}$,
then $\Pi_l$ has
three consecutive vertices with label $i$. They determine
 two consecutive edges of $\Pi_l$, which are called proper pairs.
\item[(iv)] If $i\neq j$ are
the labels of the end points of an edge $e$ of some $\Pi_l$,
then there exists exactly one more edge
among the all together $3m$ edges, whose
endpoints are labeled $i$ and $j$. This edge is called
the proper pair of $e$.
\item[(v)] It never occurs that for some
$i\neq j$, the two neighbours of a vertex
with label $i$ are both labeled $j$.
\end{description}

Given a proper labeling,
we give a orientation to each edge.
For an edge $e$ of some $\Pi_l$, whose end
points are labeled $i$ and $j$,
we orient the edge $e$ according to the natural ordering
of the labels of its vertices
if $i\neq j$, and according to the positive orientation
of $\Pi_l$ if $i=j$. We say that the orientation
of $e$ is positive if it coincides
with the orientation induced by the positive
orientation of $\Pi_l$, and the orientation
of $e$ is negative otherwise.
We say that the proper labeling
is oriented, if for any proper pair $e$ and $f$ of edges,
the orientations of $e$ and $f$ are opposite.

For $k\geq 7$, let us assume
that $\Pi_1,\ldots,\Pi_n$ are regular $k$-gons
with angle $2\pi/k$, and we have
have a proper labeling of the vertices.
Then identifying the proper pairs of edges
according their orientations,
we obtain a
(connected) compact hyperbolic surface $X$,
which is oriented if the proper labeling is oriented.

From a given proper labeling of
$n$ polygons, we will construct a proper labeling
of $n$ polygons with a higher number of vertices
using the operations (a), (b) and (c) below.

A polygonal path $a_1\ldots a_m$, $m\geq 2$,
is the union of the segments with end points
$a_i$ and $a_{i+1}$, $i=1,\ldots,m-1$, where
$a_1,\ldots,a_m$ are different.
 When we replace
a polygonal path $a_1\ldots a_m$ contained
in the boundary of a convex polygon $\Pi$
with the polygonal
path $b_1\ldots b_l$ with $a_1=b_1$
and $a_m=b_l$, we always do it
in a way  to obtain a new convex polygon $\Pi'$,
where $b_1,\ldots,b_l$ are vertices of $\Pi'$.

Let us assume that we have a proper labeling
of vertices of the convex polygons
$\Pi_1,\ldots,\Pi_n$ using labels from $\{1,\ldots,m\}$.
When replacing a polygonal path with a new one,
we only list the corresponding labels.
\begin{description}
\item[(a) at $w$ ] Assume that for
$\{x,y,z,w\}\subset \{1,\ldots,m\}$,
the union of the boundaries of $\Pi_1,\ldots,\Pi_n$
contains polygonal paths with labels
$xwy$, $ywz$ and $zyx$. Then
for $\alpha=m+1$ and $\beta=m+2$, we replace (see Figure~1)
\begin{eqnarray*}
xwy &\mbox{ by }&x\alpha w\beta y\\
ywz &\mbox{ by }&y\beta \alpha w z\\
zwx &\mbox{ by }&z w \beta \alpha x.
\end{eqnarray*}
\begin{figure}
\begin{center}
\includegraphics[width=15em]{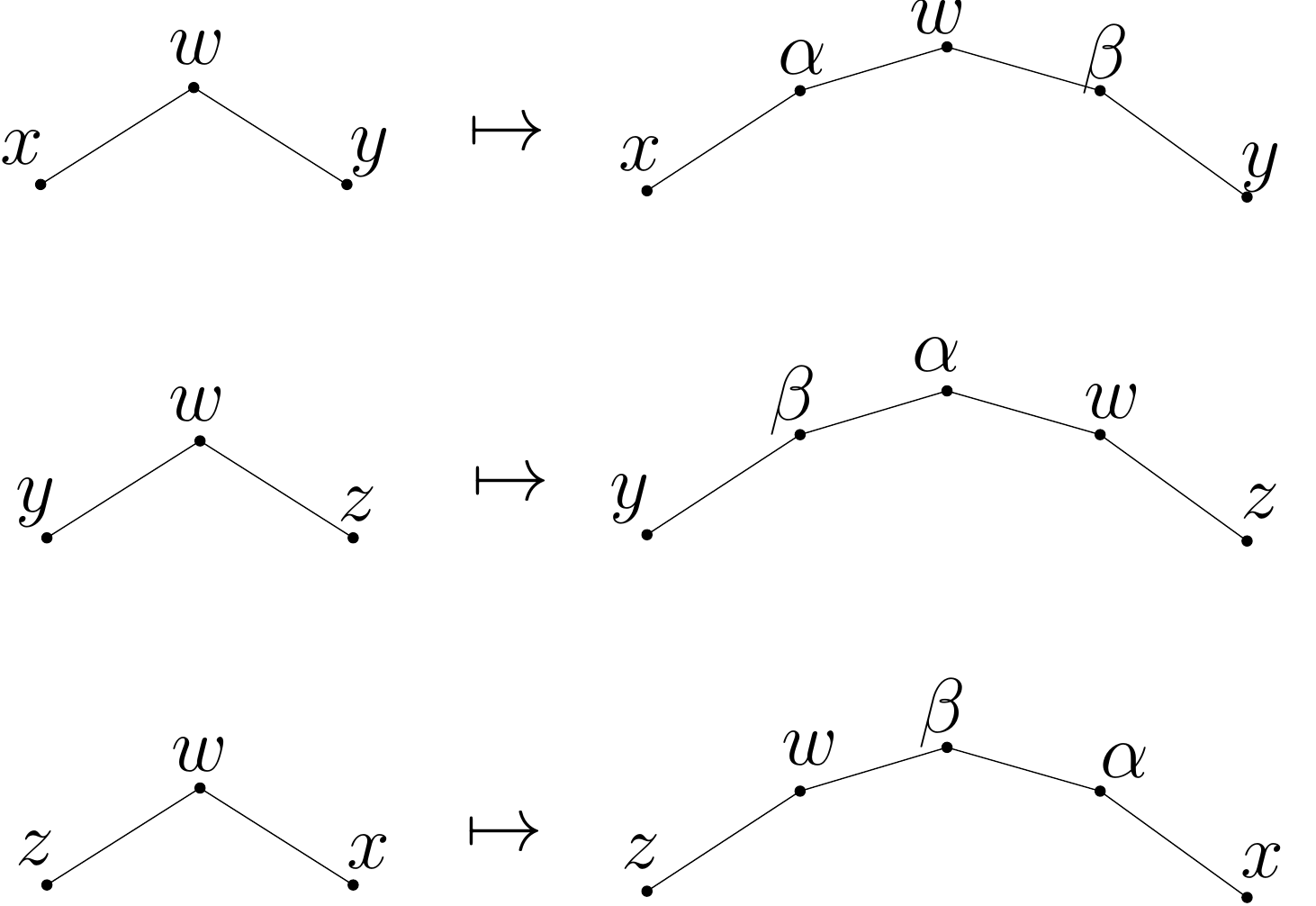}
\end{center}
\caption{ }
\end{figure}
\item[(b) at $w$ ] Assume that for
$\{x,y,z,w\}\subset \{1,\ldots,m\}$,
the union of the boundaries of $\Pi_1,\ldots,\Pi_n$
contains polygonal paths with labels
$xwy$, $ywz$ and $zyx$. Then
for $\alpha=m+1$, $\beta=m+2$,
$\gamma=m+3$ and $\delta=m+4$, we replace (see Figure~2)
\begin{eqnarray*}
xwy &\mbox{ by }&x\alpha\gamma w\delta\beta y\\
ywz &\mbox{ by }&y\beta\gamma\delta \alpha w z\\
zwx &\mbox{ by }&z w\gamma \beta\delta \alpha x.
\end{eqnarray*}
\begin{figure}
\begin{center}
\includegraphics[width=20em]{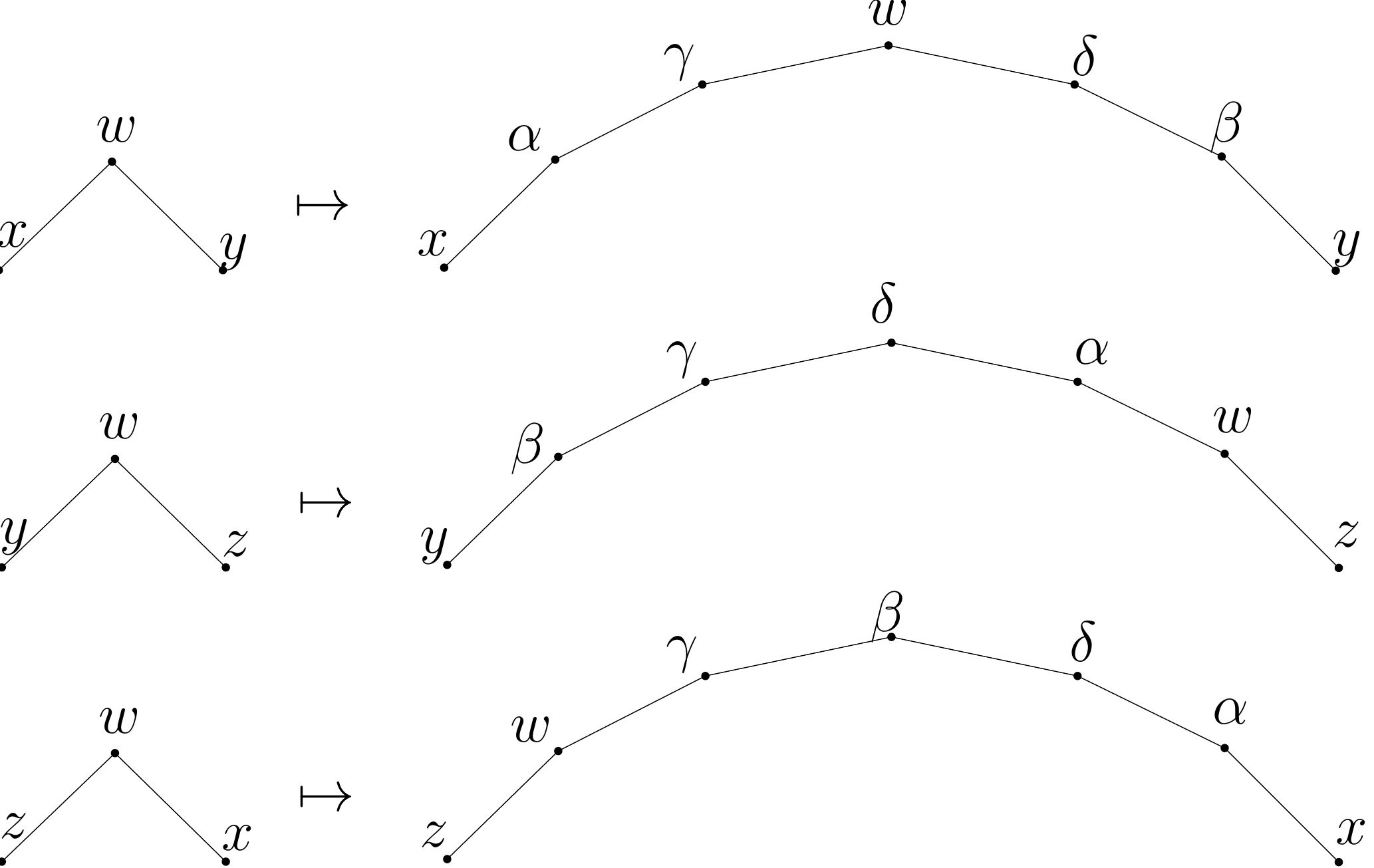}
\end{center}
\caption{ }
\end{figure}

\item[(c) at $x,y$ ] Assume
the endpoints of the proper
pair of edges $e$ and $f$ are labeled $x\neq y$.
For $\alpha=m+1$ and $\beta=m+2$, we replace (see Figure~3)
\begin{eqnarray*}
e  &\mbox{ by the polygonal path }&x\alpha \beta\beta\beta \alpha y\\
f &\mbox{ by the polygonal path }& x \alpha y.
\end{eqnarray*}
\begin{figure}
\begin{center}
\includegraphics[width=20em]{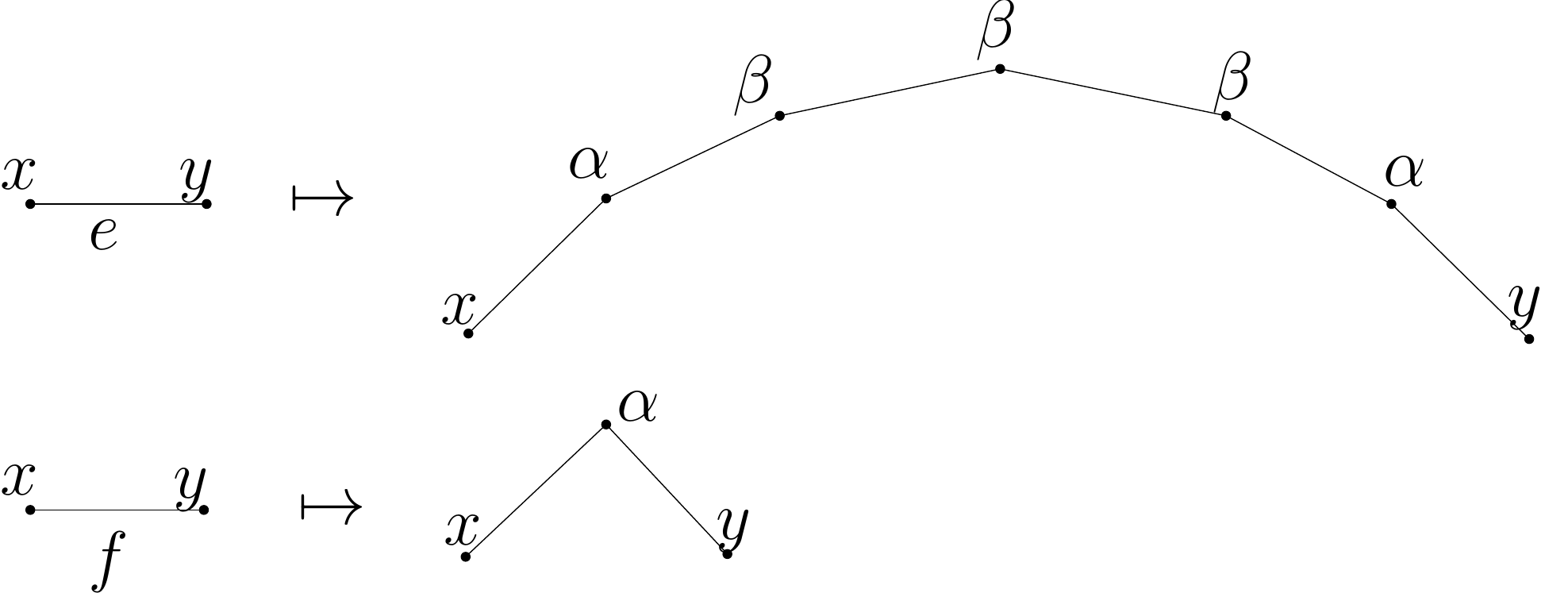}
\end{center}
\caption{ }
\end{figure}

\end{description}

We observe that all operations (a), (b) and (c)
produce a new proper labeling.
In addition,
if the proper labeling
of the vertices of $\Pi_1,\ldots,\Pi_n$
is oriented, then the new proper labeling
constructed in (b) is also oriented.

\section{Proof of Theorem~\ref{main}}

For $k\geq 7$, let $N(k)$ be the smallest
positive integer such that $N(k)\cdot k$ is divisible by six.
It follows by Lemma~\ref{nk6} that
if a compact hyperbolic surface  $X$
is tiled by embedded regular
triangles of angle $2\pi/k$, $k\geq 7$, then
\begin{equation}
\label{lN}
\chi(X)=l\cdot N(k)\cdot\left(1-\frac{k}6\right)
\mbox{ \ for a positive integer $l$}.
\end{equation}
In addition, let $X_k$ ($\widetilde{X}_k$)
be a hyperbolic surface (oriented hyperbolic surface)
of smallest area
that can be tiled by embedded regular triangles of angle $2\pi/k$.
Theorem~\ref{main} is equivalent proving that
\begin{equation}
\label{chiXk}
\mbox{there exists a proper labeling
of $N(k)$ $k$-gons.}
\end{equation}
It follows by Lemma~\ref{nonoriented} and (\ref{lN}) that
\begin{equation}
\label{chiXkori}
\mbox{if $\chi(X_k)=N(k)(1-\frac{k}6)$
and $\chi(X_k)$ is odd then $\chi(\widetilde{X}_k)=2\chi(X_k)$}.
\end{equation}
We construct $X_k$ and $\widetilde{X}_k$, depending on
the remainder of $k$ modulo $12$, and providing
a proper labeling of the vertices of $N(k)$
$k$-gons. We only provide the $N(k)$ lists
of labels according to
positive orientation of the $k$-gons.\\

\noindent{\bf Case 1 } $k\equiv 1,5,7,11\,{\rm mod}\,12$

In these cases, we have  $N(k)=6$. If $k=7$, then the following
is a proper labeling (see Figure~4):
\begin{figure}
\begin{center}
\includegraphics[width=11em]{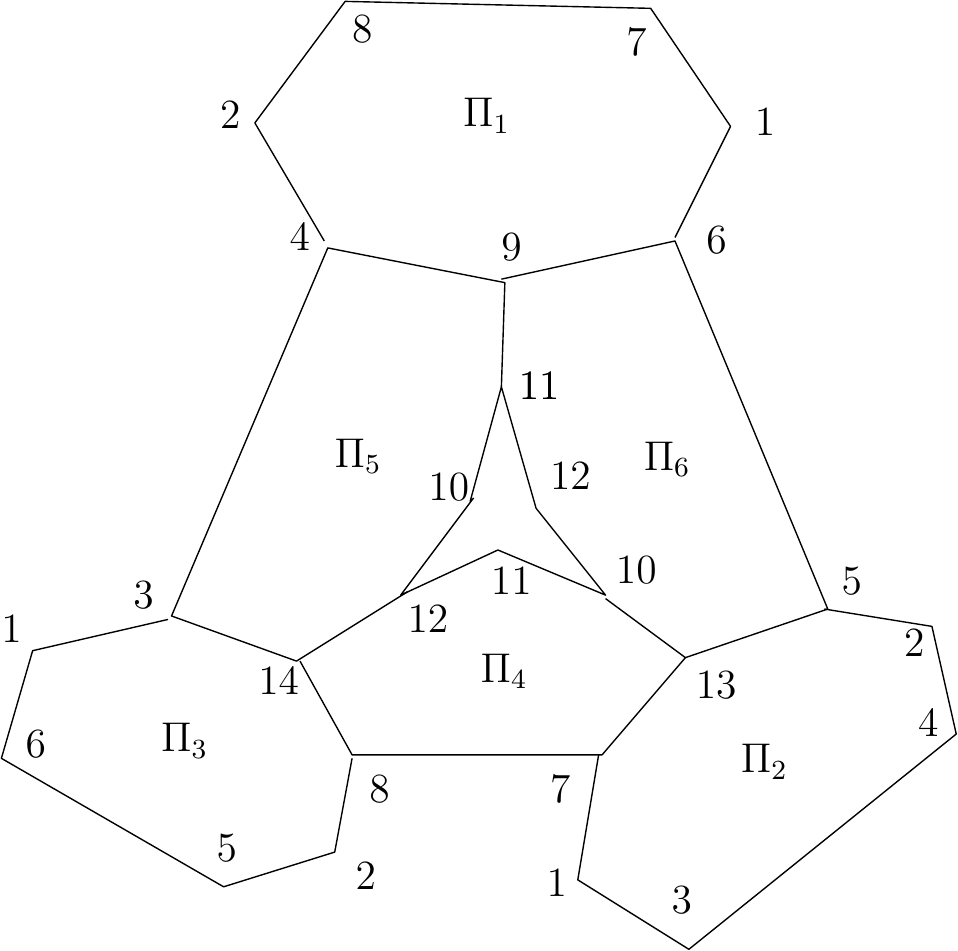}
\end{center}
\caption{ }
\end{figure}
\begin{eqnarray*}
\Pi_1: && 1,6,9,4,2,8,7\\
\Pi_2: && 1,7,13,5,2,4,3\\
\Pi_3: && 1,3,14,8,2,5,6\\
\Pi_4: && 7,8,14,12,11,10,13\\
\Pi_5: && 3,4,9,11,10,12,14\\
\Pi_6: && 9,6,5,13,10,12,11
\end{eqnarray*}
To prove (\ref{chiXk}) for $k=11$, we apply operation
(b) at $1$ (to alter $\Pi_1$, $\Pi_2$ and $\Pi_3$),
and then again (b) at $10$ (to alter $\Pi_4$, $\Pi_5$ and $\Pi_6$).
This way we obtain a proper labeling of some six
$11$-gons. After this,
we apply (a) at $1$ and $10$ to obtain a proper labeling of some six
$13$-gons. Continuing this way, applying alternately
operation (b) at  $1$ and $10$, and operation (a) at  $1$ and $10$,
we obtain a proper labeling  of some six
$k$-gons for any $k$ in Case 1.

For these $k$, we have $\chi(X_k)=6-k$ by Lemma~\ref{nk6},
which is odd.
Therefore, $\chi(\widetilde{X}_k)=2\chi(X_k)$
by (\ref{chiXkori}).\\

\noindent{\bf Case 2 } $k\equiv 2,10\,{\rm mod}\,12$

In these cases, we have  $N(k)=3$. If $k=10$, then the following
is an oriented proper labeling (see Figure~5):
\begin{figure}
\begin{center}
\includegraphics[width=11em]{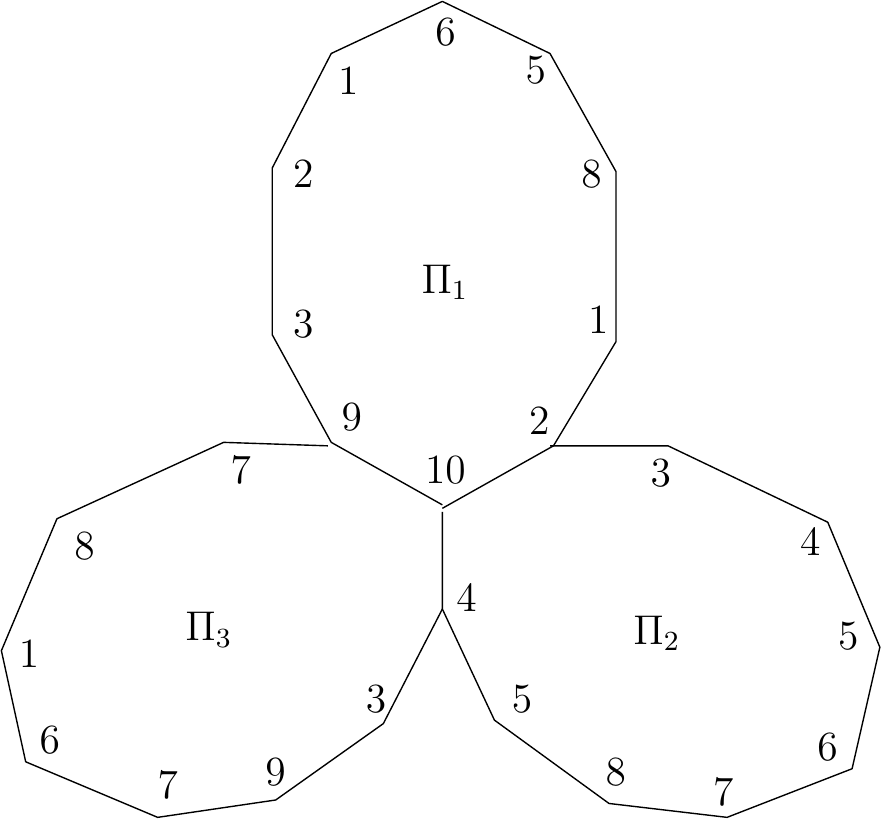}
\end{center}
\caption{ }
\end{figure}
\begin{eqnarray*}
\Pi_1: && 10,9,3,2,1,6,5,8,1,2\\
\Pi_2: && 10,2,3,4,5,6,7,8,5,4\\
\Pi_3: && 10,4,3,9,7,6,1,8,7,9
\end{eqnarray*}
To prove (\ref{chiXk}) for $k=14$, we apply operation
(b) at $10$. To prove (\ref{chiXk}) for $k=22$,
we apply operation
(b) at $10$ twice. In particular, using operation
(b) at $10$, we obtain a proper labeling  of some three
$k$-gons for any $k$ in Case 2.

Since each one of the proper labelings presented
in Case 2 is oriented, the $X_k$ constructed this way
is also oriented. \\

\noindent{\bf Case 3 } $k\equiv 3,9\,{\rm mod}\,12$

In these cases, we have  $N(k)=2$. If $k=9$, then the following
is a proper labeling (see Figure~6):
\begin{figure}
\begin{center}
\includegraphics[width=11em]{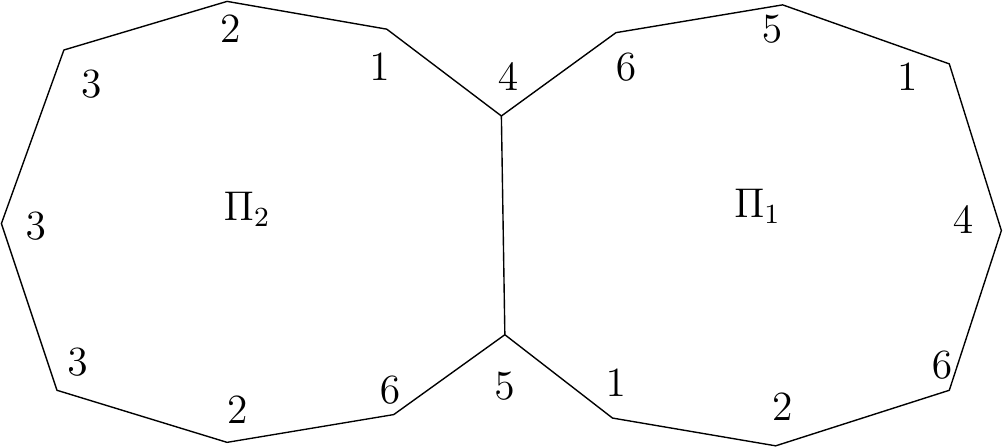}
\end{center}
\caption{ }
\end{figure}
\begin{eqnarray*}
\Pi_1: && 4,6,5,1,4,6,2,1,5\\
\Pi_2: && 4,5,6,2,3,3,3,2,1.
\end{eqnarray*}
To prove (\ref{chiXk}) for $k=15$, we apply operation
(c) at $4,5$. In particular, using operation
(b) at $4,5$, we obtain a proper labeling  of some pair of
$k$-gons for any $k$ in Case 3.

For these $k$, we have $\chi(X_k)=2-\frac{k}3$ by Lemma~\ref{nk6},
which is odd as $\frac{k}3$ is odd.
Therefore, $\chi(\widetilde{X}_k)=2\chi(X_k)$
by (\ref{chiXkori}).\\

\noindent{\bf Case 4 } $k\equiv 4,8\,{\rm mod}\,12$

In these cases, we have  $N(k)=3$. If $k=8$, then the following
is a proper labeling (see Figure~7):
\begin{figure}
\begin{center}
\includegraphics[width=11em]{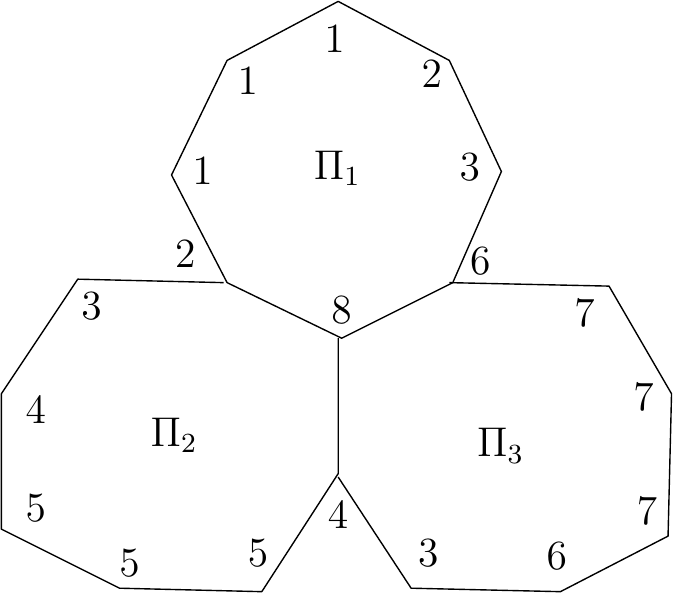}
\end{center}
\caption{ }
\end{figure}
\begin{eqnarray*}
\Pi_1: && 8,2,1,1,1,2,3,6\\
\Pi_2: && 8,4,5,5,5,4,3,2\\
\Pi_3: && 8,6,7,7,7,6,3,4.
\end{eqnarray*}
To prove (\ref{chiXk}) for $k=16$, we apply operation
(b) at $8$ twice. To prove (\ref{chiXk}) for $k=20$,
we apply operation
(b) at $8$. In particular, using operation
(b) at $8$, we obtain a proper labeling  of some three
$k$-gons for any $k$ in Case 4.

For these $k$, we have $\chi(X_k)=3-\frac{k}2$ by Lemma~\ref{nk6},
which is odd as $\frac{k}2$ is even.
Therefore, $\chi(\widetilde{X}_k)=2\chi(X_k)$
by (\ref{chiXkori}).\\

\noindent{\bf Case 5 } $k\equiv 6\,{\rm mod}\,12$

In this case, we have $N(k)=1$. If $k=18$, then the following
is an oriented proper labeling (see Figure~8):
\begin{figure}
\begin{center}
\includegraphics[width=11em]{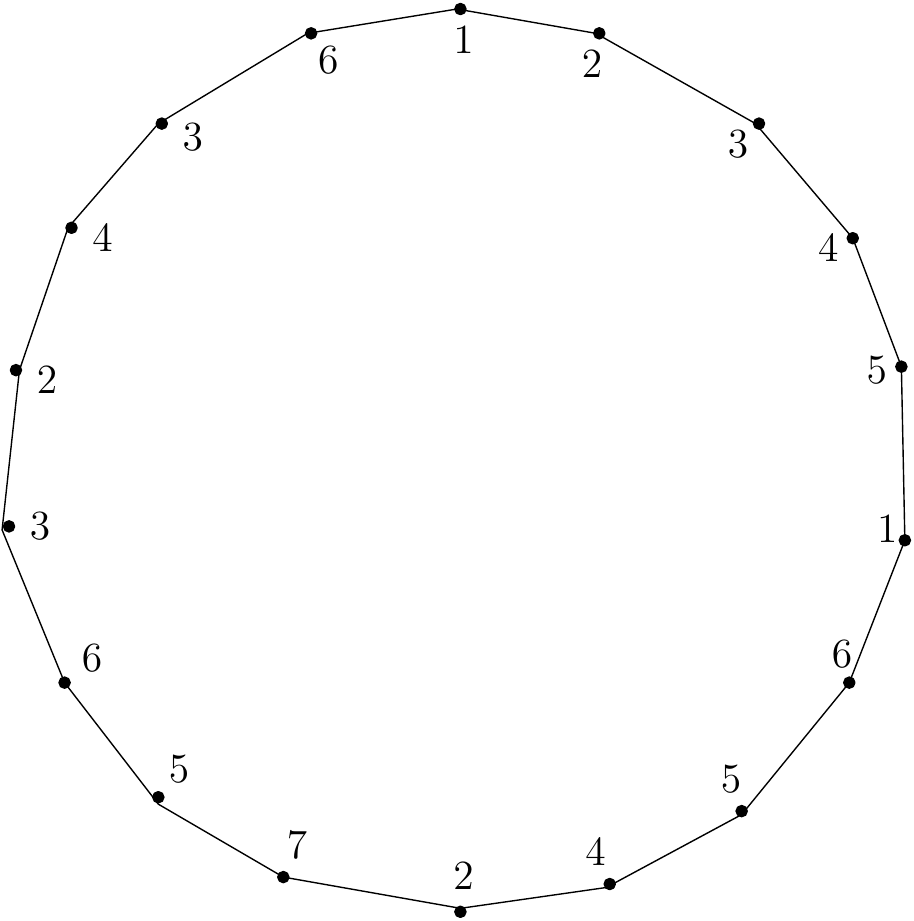}
\end{center}
\caption{ }
\end{figure}
\begin{eqnarray*}
\Pi_1: && 1,2,3,4,5,1,6,5,4,2,7,5,6,3,2,4,3,6.
\end{eqnarray*}
To prove (\ref{chiXk}) for $k=30$, we apply operation
(b) at $1$ three times. In particular, using operation
(b) at $1$, we obtain an oriented proper labeling  of some
$k$-gon for any $k$ in Case 5.
Therefore, $\chi(\widetilde{X}_k)=\chi(X_k)=1-\frac{k}6$ in this case.\\

\noindent{\bf Case 6 } $k\equiv 0\,{\rm mod}\,12$

In this case, we have $N(k)=1$. If $k=12$, then the following
is a proper labeling (see Figure~9):
\begin{figure}
\begin{center}
\includegraphics[width=11em]{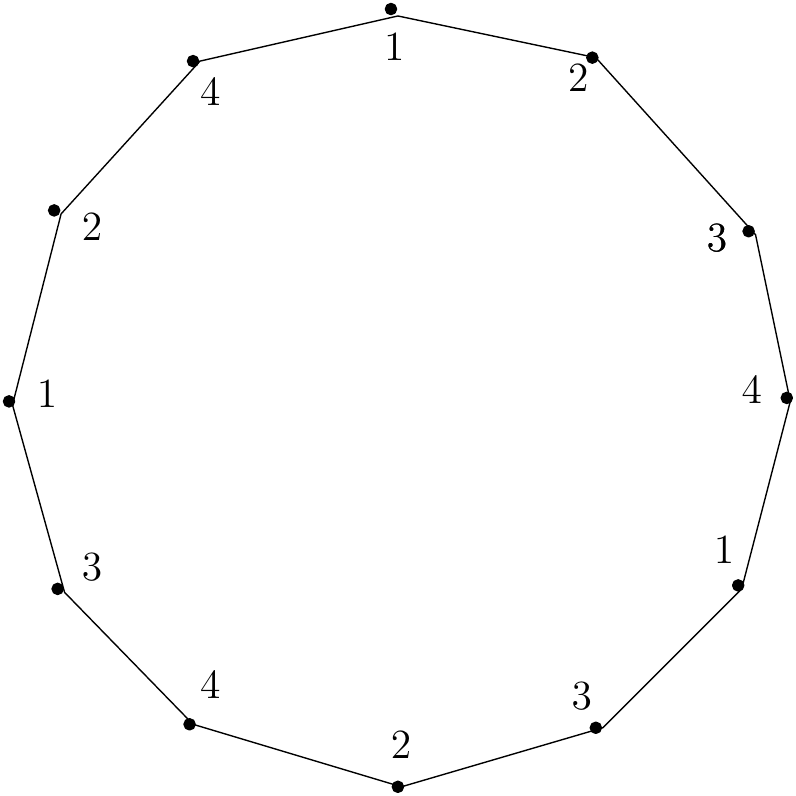}
\end{center}
\caption{ }
\end{figure}
\begin{eqnarray*}
\Pi_1: && 1,2,3,4,1,3,2,4,3,1,2,4.
\end{eqnarray*}
To prove (\ref{chiXk}) for $k=24$, we apply operation
(b) at $1$ three times. In particular, using operation
(b) at $1$, we obtain a proper labeling  of some
$k$-gon for any $k$ in Case 6.

For these $k$, we have $\chi(X_k)=1-\frac{k}6$ by Lemma~\ref{nk6},
which is odd as $\frac{k}6$ is even.
Therefore, $\chi(\widetilde{X}_k)=2\chi(X_k)$
by (\ref{chiXkori}). Q.E.D.\\

\noindent {\bf Remark } If $k\geq 7$ is divisible by six, then there
is a graph theoretic way to produce a proper labeling of a
$k$-gon. Let $G$ be a connected graph on $\{1,\ldots,m\}$ for $m=k/3$
such that each vertex is of degree three. 
We call a closed path on the edges of $G$ 
a double Hamiltonian path if each edge is traveled exactly twice,
and consecutive edges in the path are always different. 
Then
a proper labeling of the vertices of the $k$-gon 
where no pair of consecutive vertices have the same label, is
equivalent to a double Hamiltonian path. The proper labeling
is oriented if and only if each edge
is traveled in both directions in the double Hamiltonian path.

%\noindent {\bf Acknowledgement: }

\noindent Christophe Bavard, {\it Christophe.Bavard@math.u-bordeaux.fr }\\
Institut de Math\'ematiques, Universit\'e Bordeaux 1\\

\noindent K\'aroly J. B\"or\"oczky, {\it carlos@renyi.hu}\\
Alfr\'ed R\'enyi Institute of Mathematics, and\\
Universitat Polit\`ecnica de Catalunya, Barcelona Tech, and\\
 Department of Geometry, Roland E\"otv\"os University \\

\noindent Borb\'ala Farkas, {\it farkas\_borka@yahoo.co.uk}\\
 Department of Geometry, Technical University of Budapest\\

\noindent Istv\'an Prok, {\it prok@math.bme.hu}\\
 Department of Geometry, Technical University of Budapest\\

\noindent Lluis Vena, {\it lluis.vena@utoronto.ca}\\
Department of Mathematics, University of Toronto\\

\noindent Gergely Wintsche, {\it wgerg@ludens.elte.hu}\\
Department of Didactics, Roland E\"otv\"os University \\


\begin{thebibliography}{10}

\bibitem{AVS93}
D.V. Alekseeskij, E.B. Vinberg, A.S. Solodovnikov:
 Geometry of spaces of constant curvature.
Encycl. Math. Sci. 29, 
Springer--Verlag, 1993, 1-138.

\bibitem{Bav96}
C. Bavard:
Disques extr\'emaux et surfaces modulaires.
Ann. Fac. Sci. Toulouse Math.,  5 (1996), 191-202.

\bibitem{EEK82a}
A.L. Edmonds, J.H. Ewing, R.S. Kulkarni:
Regular tessellations of surfaces and (p, q, 2)-triangle groups.
Ann. Math., 116 (1982), 113-132.


\bibitem{EEK82i}
A.L. Edmonds, J.H. Ewing, R.S. Kulkarni:
Torsion Free Subgroups of Fuchsian Groups and Tessalations of Surfaces.
Inventiones Mathematicae, 69 (1982), 331-346.

\bibitem{FTK93}
G. Fejes T\'oth, W. Kuperberg:  
Packing and covering with convex sets,   
In: Handbook of Convex Geometry, 
P. Gruber and J. Wills (eds.), North-Holland 1993, 799-860. 

\bibitem{FTL64}
L. Fejes T\'oth:
Regular Figures.
Pergamon Press, 1964.

\bibitem{Rat94}
J.G. Ratcliffe:
Foundations of hyperbolic manifolds.
Springer, 1994.

\end{thebibliography}
\end{document}